\documentclass[a4paper,12pt,twoside]{article}

\usepackage{amsfonts,amsthm,amsmath,amsfonts,latexsym,amssymb}

\usepackage[ansinew]{inputenc}
\usepackage{graphicx,upref,calligra}
\usepackage[svgnames]{xcolor}
\usepackage{mathrsfs}
\usepackage[T1]{fontenc}
\usepackage{enumerate}
\usepackage{esint}
\usepackage{authblk}

\usepackage{fancyhdr}
\usepackage{lastpage}

\newtheorem{fed}{Definition}[section]

\newtheorem{teo}[fed]{Theorem}
\newtheorem*{teo*}{Theorem}
\newtheorem{lem}[fed]{Lemma}
\newtheorem{cor}[fed]{Corollary}
\newtheorem{pro}[fed]{Proposition}
\theoremstyle{definition}
\newtheorem*{rem}{Remark}

\newtheorem*{exa}{Example}

\definecolor{titleblue}{rgb}{0.13,0.49,0.69}
\definecolor{mylred}{rgb}{0.85,0.24,0.2}
\definecolor{myyellow}{rgb}{0.42,0.24,0.52}
\definecolor{mygreen}{rgb}{0.12,0.5,0.29}
\definecolor{myred}{rgb}{0.74,0.13,0.13}
\definecolor{mylblue}{rgb}{0.2,0.75,1}
\definecolor{mylgreen}{rgb}{0.68,0.98,0.6}
\definecolor{mylyellow}{rgb}{0.86,0.85,0.55}
\definecolor{myllyellow}{rgb}{0.87,0.86,0.56}
\definecolor{naranja}{RGB}{249,153,96}
\definecolor{sidebardarkcolor}{rgb}{0.21,0.31,0.40}
\definecolor{sidebarlightcolor}{rgb}{0.7,0.77,0.836}

\def\EOE{\hfill $\blacktriangle$}

\def\bdem{\begin{proof}}
\def\edem{\end{proof}}
\def\ds{\displaystyle}

\def\bm{\left(\begin{array}}
\def\em{\end{array}\right)}
\def\ben{\begin{enumerate}}
\def\een{\end{enumerate}}
\def\barr{\begin{array}}
\def\earr{\end{array}}

\def\eps{\varepsilon}
\def\fii{\varphi }

\def\w{\omega}
\def\W{\Omega}

\def\N{\mathbb{N}}
\def\Z{\mathbb{Z}}
\def\R{\mathbb{R}}
\def\C{\mathbb{C}}

\def\T{\mathbb{T}}

\def\cB{\mathcal{B}}

\def\cP{\mathcal{P}}

\newcommand{\peso}[1]{ \quad \text{ \rm  #1 } \quad }

\newcommand{\sub}[2]{{#1}_{\mbox{\tiny{${#2}$}}}}

\newcommand{\barc}[1]{{\beta}_{#1}}
\newcommand{\barza}[1]{\sub{\beta}{#1}}

\newcommand{\mat}{\mathcal{M}_n(\mathbb{C})}

\newcommand{\matsa}{\mathcal{H}(n)}

\newcommand{\matpos}{\mat^+}

\newcommand{\conv}{\xrightarrow[n\rightarrow\infty]{}}

\numberwithin{equation}{section}

\newcommand{\argmin}{\arg\!\min}

\begin{document}

\title{\textbf{Ergodic theorem in Hadamard spaces in terms of inductive means}}
\author[1,2]{Jorge Antezana}
\author[1, 2]{Eduardo Ghiglioni}
\author[1, 2]{Demetrio Stojanoff}
\affil[1]{Departamento de de Matem\'atica, FCE-UNLP, Calles 50 y 115, 
	(1900) La Plata, Argentina.}
\affil[2]{ Instituto Argentino de Matem\'atica, `Alberto P. Calder\'on', CONICET, Saavedra 15 3er. piso,
	(1083) Buenos Aires, Argentina.}

\date{}

	\maketitle
	
	\begin{abstract}
		Let $(G,+)$ be a compact, abelian, and metrizable topological group. In this group we take $g\in G$ such that the corresponding automorphism $\tau_g$ is ergodic. The main result of this paper is a new ergodic theorem for functions in $L^1(G,M)$, where $M$ is a Hadamard space. The novelty of our result is that we use inductive means to average the elements of the orbit $\{\tau_g^n(h)\}_{n\in\mathbb{N}}$. The advantage of inductive means is that they can be explicitly computed in many important examples. The proof of the ergodic theorem is done firstly for continuous functions, and then it is extended to $L^1$ functions. The extension is based in a new
		construction of mollifiers in Hadamard spaces. This construction has the advantage that it only uses the metric structure and the existence of barycenters, and do not require the existence of an underlyng vector space. For this reason, it can be used in any Hadamard space, in contrast with those results that need to use the tangent space or some chart to define the mollifier.
	\end{abstract}

	\section{Introduction}
	
	\subsection{Motivation}
	The initial motivation for this work was the study of ergodic theorems involving the barycenters in the space of $n\times n$ strictly positive matrices $\matpos$. More precisely, we pursued an ergodic type theorem in terms of the so called inductive means, which in many important cases can be explicitly computed. 
	
	\medskip
	
	Recall that the set $\matpos$ is an open cone in the real vector space of selfadjoint matrices $\matsa$. In particular, it is a differential manifold and the tangent spaces can be
	indentified for simplicity with $\matsa$. The manifold $\matpos$ can be endowed with a natural Riemannian structure. With respect to this metric structure,  if  $\alpha:[a,b]\to\matpos$ is a piecewise smooth path, its length is defined by 
	$$
	L(\alpha)=\int_a^b \|\alpha^{-1/2}(t)\alpha'(t)\alpha^{-1/2}(t)\|_2 \,dt,
	$$
	where $\|\cdot\|_2$ denotes the Frobenius or Hilbert-Schmidt norm.
	%
	In this way, $\matpos$ becomes a Riemannian manifold with non-positive curvature. As usual, a distance $\delta$ can be defined by 
	$$
	\delta(A,B)=\inf \{L(\alpha):\ \mbox{$\alpha$ is a piecewise smooth path connecting $A$ with $B$}\}.
	$$
	The infimum is actually a minimum, and the geodesic connecting two positive matrices $A$ and $B$ has the following simple expression
	$$
	\gamma_{AB}(t)=A^{1/2}(A^{-1/2}BA^{-1/2})^{t}A^{1/2}\,.
	$$
	It is usual in matrix analysis to use the notation $A\#_t B$ instead of $\gamma_{AB}(t)$. The midpoint $A\#_{\frac{1}{2}} B$ is called \textbf{geometric mean} or \textbf{barycenter} between $A$ and $B$,  and it also admits the following variational characterization
	$$
	A\#_{\frac12} B=\argmin_{C\in \matpos} \ \Big(\delta^2(A,C)+\delta^2(B,C)\,\Big).
	$$
	There is no reason to restrict our attention to only two matrices. The notion of geometric mean can be generalized for more than two matrices in the obvious way
	$$
	\Gamma(A_1,\ldots,A_n):=\argmin_{C\in \matpos} \ \Big(\sum_{j=1}^n\delta^2(A_j,C)\,\Big).
	$$
	The solution of this least square minimization problem exists and is unique because of the convexity properties of the distance $\delta(\cdot,\cdot)$. The geometric means naturally appear in many applied problems. For instance, they appear in the study of radar signals. In these problems, each signal is detected by more than one sensor. The information of each sensor is codified in a covariance kernel, and these kernels have to be averaged to get the final output. It turns out that the best way to average the kernels is not the standard arithmetic mean, but the geometric mean introduced above (see \cite{Barba} and the references therein for more details). Another typical application of the geometric means is in problems related with the gradient or Newton like optimization methods (see \cite{Bi-Ia},\cite{Maher}).
	
	\medskip
	
	The usual problem dealing with geometric means is that the geometric mean of three or more matrices does not have in general a closed formula.  For this reason, there has been an intensive research with the aim to find good ways to approximate the geometric mean of more than two matrices (\cite{[Bhatia2]}, \cite{[LyL]}, \cite{LimPal}). In \cite{[Holbrook]}, Holbrook proved that they can be approximated by the so called \textbf{inductive means}. To motivate the definition of inductive means, note that given a sequence $\{a_n\}_{n\in\N}$ of complex numbers
	\begin{align*}
	\frac{a_1+a_2+a_3}{3}=&\ds\frac{2}{3}\left(\frac{a_1+a_2}{2}\right)+\frac13\, a_3\\
	\vdots&\\
	\frac{a_1+\ldots+a_n}{n}=&\ds\frac{n-1}{n}\left(\frac{a_1+\ldots+a_{n-1}}{n-1}\right)+\frac1n\, a_n.\\
	\end{align*}
	Let $\gamma_{a,b}(t)=t\,b+(1-t)a$, and for a moment allow us to use the notation $a\,\sharp_t\, b=\gamma_{a,b}(t)$.  Then
	\begin{align*}
	\frac{a_1+a_2+a_3}{3}=&(a_1\,\sharp_{\frac12}\, a_2)\,\sharp_{\frac13}\,a_3\\
	\frac{a_1+a_2+a_3+a_4}{4}=&((a_1\,\sharp_{\frac12}\, a_2)\,\sharp_{\frac13}\,a_3)\,\sharp_{\frac14}\,a_4.
	\end{align*}
	and so on and so forth. The segments are the geodesics in the euclidean space. Thus, in our setting, we can replace the segments by the geodesic associated to the Riemannian structure. This is the idea that leads to the definition of the inductive means. Given a sequence of strictly positive matrices $\{A_n\}_{n\in\N}$, the inductive means are defined as follows:
	\begin{align*}
	\Gamma_1(A) & = A_1\\
	\Gamma_{n}(A) & = \Gamma_{n-1}(A)  \#_{\frac{1}{n}} A_n \ \ \ \ \ (n \geq 2).
	\end{align*}
	Now, consider $d$ positive matrices $A_0,\ldots,A_{d-1}$, and define the function $F:\Z_d\to\matpos$ by $F(\overline{k})=A_k$, where $\Z_d$ denotes the abelian group of integers mod $d$. Then, if we define the periodic sequence $A=\{F(\overline{n})\}_{n\in\N}$, then Holbrook proved 
	$$
	\lim_{n\to\infty} \Gamma_{n}(A) =\Gamma(A_0,\ldots,A_{d-1}).
	$$

	\subsection{About our results}
	
	Let $(G,+)$ be a compact topological group, endowed with a Haar measure $m$, and let $\tau:G\to G$ be an ergodic map. The classical Birkhoff ergodic theorem says that, given $f\in L^1(G)$, then for $m$-almost every $h\in G$
	\begin{equation}\label{Birky}
	\frac{1}{n}\sum_{k=0}^{n-1} f(\tau^k(\w))\conv \int_{G} f(\w) dm(\w).
	\end{equation}
Holbrook's theorem can be seen as an ergodic theorem for functions defined in $\Z^p$ and taking values in $\matpos$. 
Indeed, let the inductive means play the role of the averages at the left hand side, and let the joint geometric mean $\Gamma(A_0,\ldots,A_{d-1})$ play the role of the integral of $F$ with respect to the Haar measure in $\Z_d$. 

\medskip

This interpretation of Holbrook result suggests that we can get computable approximations of the geometric means of several matrices in terms of ergodic averages. Our main goal is to prove such ergodic theorems in terms of the inductive means in a much more general context. More precisely, consider a dynamical system $(G,\tau)$,  where $G$ is a compact group and $\tau$ is ergodic with respect to the Haar measure $m$. Note that the group has to be abelian because the orbits $\{\tau^{n}g\}_{n\in\N}$ are dense in $G$. In this setting recall that given a function $F:G\to\matpos$, we say that $F\in L^p(G,\matpos)$ if
	$$
	\int_G \delta^p(F(g),B))\,dm(g)<\infty,
	$$
	where $B$ is any positive matrix. By the triangular inequality, this definition does not depend on the choice of $B$. Following Sturm \cite{[sturm]}, the barycenter of $F\in L^1(G,\matpos)$ is defined as
	$$
	\barza{F}:=\argmin_{C\in\matpos}\int_G \delta^2(F(g),C))-\delta^2(F(g),B))\,dm(g).
	$$ 
	As before, this definition  does not depend on $B$. Note that if $F\in L^2(G,\matpos)$ 
	$$
	\barza{F}:=\argmin_{C\in\matpos}\int_\T \delta^2(F(g),C))dm(g)
	$$
	and we get a natural generalization of the geometric means defined before. 
	In this paper we  prove that, given $A\in L^1(G, \matpos)$, for almost every $g\in G$
	\begin{equation}\label{eq intro ergodic}
	\lim_{n\to\infty} \Gamma_{n}\big(A(g),A(\tau(g)),\ldots,A(\tau^{n-1}(g))\big)=\barza{A}.
	\end{equation}
	
	Moreover, we prove this result not only for functions taking values in $\matpos$, but also in any Hadamard space $M$(see subsection \ref{hadamard} for a precise definition).  In \cite{[sturm]},  Sturm developed a theory of barycenters of probability measures for Hadamard spaces (see Subsection \ref{Baricentro} for some definitions and basic results). Endowed with this barycenter, Hadamard spaces play an important role in the theory of integrations (random variables, expectations and variances), law of large numbers, ergodic theory, Jensen's inequality (see \cite{Bochi},  \cite{Sahib}, \cite{LL4},  \cite{Navas}, and \cite{[sturm]}), stochastic generalization of Lipschitz retractions and extension problems of Lipschitz and H\"older maps (see \cite{LeeNaor}, \cite{Mendel}, and \cite{Ohta}) and optimal transport theory on Riemannian manifolds (see \cite{Pass1}, and \cite{Pass2}), etc. 
	
	\medskip
	
	The generalization from $\matpos$, or even from more general Riemannian manifolds with non-positive curvature, to general Hadamard spaces is not straightforward. One of the reasons is that a general Hadamard space does not necessarily have an underlying finite dimensional vector space, as in the case of manifolds. This makes some steps of the proof much more involved, and lead us to a new definition of mollifiers that only uses the metric structure  (see Subsection \ref{aproximando}).

%

%
	
	\medskip

	\subsection{Comments on related works}

In  \cite{Austin} Austin proved a very general ergodic theorem for Hadamard spaces, which in our setting says that, given $A\in L^2(G, M)$, for almost every $g\in G$  it holds that
	$$
	\lim_{n\to\infty} \Gamma\big(A(g),A(\tau(g)),\ldots,A(\tau^{n-1}(g))\,\big)=\barza{A}.
	$$
Later on, in \cite{Navas} Navas extended the Austin's result to functions in $A\in L^1(G, M)$, and taking values in more general metric spaces. In these results, the discrete arithmetic means of Birkhoff's theorem are replaced by the (joint) geometric mean of the $n$-tuple
	$$
	\big(A(g),A(\tau(g)),\ldots,A(\tau^{n-1}(g))\,\big).
	$$
	On the other hand, the integral in Birkhoff theorem is replaced by the barycenter $\barza{A}$. Note that in the sequence
	$$
	A(g), \,\Gamma\Big(A(g),A(\tau(g))\Big), \,\Gamma\Big(A(g),A(\tau(g),A(\tau^2(g))\Big),\,\ldots
	$$
	from the second term on, the elements of the sequence do not have a closed formula. In this direction, the inductive means are simpler and provide in many concrete instances a computable approximation sequence to the barycenter. This is an advantage of our result, but there is a price to pay for this. On one hand, we need a good control on the convexity of the metric. For this reason we work on Hadamard spaces, as in the case of Austin's result \cite{Austin}. On the other hand, the inductive means are simple averages of points. This simplicity is good from the computational point of view, but it confine our result to $\Z$-actions. The aforementioned results also hold  for more general actions.

	\subsection{Organization of the paper}
	The paper is organized as follows. Section 2 is devoted to collect some preliminaries on Hadamard spaces, as well as, barycenters and inductive means in Hadamard spaces. In section 3 we prove our main result for functions defined in a Kronecker systems and taking values in a Hadamard space. Firstly we will prove the result for continuous function (Theorem \ref{mainteo1}). In order to extend this result to $L^1$ functions (Theorem \ref{Ergodic L1}), first of all we prove in Subsection \ref{aproximando} some results related with approximation by continuous functions in general Hadamard spaces. These results are interesting by themselves, and generalize some results proved by Karcher in \cite{Karcher}. Finally, in the Subsection \ref{se termina} we complete the proof of the $L^1$ version of the ergodic theorem.

	\section{Preliminaries}
	
	\subsection{Hadamard spaces}\label{hadamard}
	
	In this section we summarize some basic facts about Hadamard spaces, also called (global) CAT(0) spaces or non-positively curved (NPC) spaces. This subject started with the works by Alexandrov \cite{Alex} and Reshetnyak  \cite{Res}. Nowadays there exists a huge bibliography on the subject. The interested reader is referred to the monographs \cite{ball}, \cite{BGS}, \cite{BH}, and \cite{[Jost]}  for more information.
	
	\medskip
	
	\begin{fed}
		A complete metric space $(M,\delta)$ is called a \textbf{Hadamard space} if it satisfies the semiparallelogram law, i.e., for each $x, y \in M$ there exists  $m \in M$ satisfying
		\begin{equation}\label{eq1}
		\delta^2(m, z) \leq \ds\frac{1}{2}\delta^2(x, z) + \ds\frac{1}{2}\delta^2(y, z) - \ds\frac{1}{4}\delta^2(x, y)
		\end{equation}
		for all $z \in M$.  The point $m$ is called (metric) \textbf{midpoint} between $x$ and $y$.
	\end{fed}
	
	Taking $z=x$ and $z=y$ in the inequality \eqref{eq1}, it is easy to conclude that $\delta(x,m)=\delta(m,y)=\frac12 \delta(x,y)$. Moreover, this inequality also implies that the midpoint is unique. The existence and uniqueness of midpoints give rise to a unique (metric) geodesic $\gamma_{a,b} : [0, 1] \rightarrow M$ connecting any given two points $a$ and $b$. Indeed, firstly define $\gamma_{a,b}(1/2)$ to be the midpoint of $a$ and $b$. Then, using an inductive argument, we define the geodesic for all dyadic rational numbers in $[0, 1]$. Finally by completeness,  it can be extended to all $t \in [0, 1]$. Throughout this paper, we will use the notation $a \#_t b$  instead of  $\gamma_{a,b}(t)$.
	
	\bigskip
	
	The inequality \eqref{eq1} also extends to arbitrary points on geodesics.
	
	\begin{pro} Let $(M,\delta)$ be a Hadamard space. Then, for all $t \in [0, 1]$ and $x, y, z \in M$,
		\begin{equation}\label{eq2}
		\delta^2(x \#_t y, z) \leq (1-t)\delta^2(x, z) + t\delta^2(y, z) - t(1-t)\delta^2(x, y).
		\end{equation}
	\end{pro}
	
	\medskip
	
	\noindent A consequence of this result, that we will use later, is:

	\begin{cor}\label{cor4}
		
		Given four points $a, a^{'}, b, b^{'} \in M$ let
		$$
		f(t) = \delta(a \#_t a^{'}, b \#_t b^{'}).
		$$
		Then $f$ is convex on $[0, 1]$; i.e.
		\begin{equation}\label{eq5}
		\delta(a \#_t a^{'}, b \#_t b^{'}) \leq (1-t)\delta(a, b) + t\delta(a^{'}, b^{'}).
		\end{equation}
		
	\end{cor}
	
	\bigskip
	
	\noindent We conclude this subsection with the so called Reshetnyak's Quadruple Comparison theorem.
	\begin{teo}\label{Res}Let $(M, \delta)$ be a Hadamard space. For all $x_1, x_2, x_3, x_4 \in M$,
		\begin{equation}\label{res}
		\delta^2(x_1, x_3) + \delta^2(x_2, x_4) \leq \delta^2(x_2, x_3) + \delta^2(x_1, x_4) + 2\delta(x_1, x_2)\delta(x_3, x_4). 
		\end{equation}	
	\end{teo}
	
	%
	%
	%
	%
	%
	
	\subsection{Barycenters}\label{Baricentro}
	
	Let $(M, \delta)$ be a Hadamard space, and let $\cB(M)$ the $\sigma$-algebra of Borel sets (i.e. the smallest $\sigma$-algebra that contains the open sets).
	Denote by  $\mathcal{P}(M)$ the set of all probability measures on $\mathcal{B}(M)$ with separable support, and for $1 \leq \theta < \infty$, let $\mathcal{P}^{\theta}(M)$ denote the set of $\mu \in \mathcal{P}(M)$ such that
	$$
	\int \delta^{\theta}(x, y)d\mu(y) < \infty,
	$$
	for some (hence all) $x \in M$. By means of $\mathcal{P}^{\infty}(M)$  we will denote the set of all measures in $\mathcal{P}(M)$ with bounded support. 
	
	\medskip
	
	\begin{pro}\label{bary}
		Let $(M, \delta)$ be a Hadamard space and fix $y \in M$. For each $\mu \in \mathcal{P}^1(M)$ there exists a unique point $\barc{\mu} \in M$ which minimizes the uniformly convex, continuous function
		$$
		z \mapsto \int_M [\delta^2(z, x) - \delta^2(y, x)]d\mu(x).
		$$
		This point is independent of $y$.
	\end{pro}
	
	Following Sturm's paper \cite{[sturm]}, the point $\barc{\mu}$ is called \textbf{barycenter} of $\mu$. 
	If $\mu \in \mathcal{P}^2(M)$ then $\barc{\mu}$ coincides with the usual Cartan's definition of barycenter:
	$$
	\argmin_{z \in M} \int_M \delta^2(z, x)d\mu(x).
	$$
	
	The following inequality satisfied by the baricenter will be very important in the sequel.
	
	\begin{pro}[Variance Inequality]\label{baryvar}
		Let $(M, \delta)$ be a Hadamard space. For any probability measure $\mu \in \mathcal{P}^1(M)$ and for all $z \in M$:
		\begin{equation}\label{semip}
		\int_{M} [\delta^2(z, x) - \delta^2(\barc{\mu}, x)]d\mu(x) \geq \delta^2(z, \barc{\mu}).
		\end{equation}
	\end{pro}

	\bigskip
	
	Now, let $(\Omega, P)$ be an arbitrary probability space and let $F: \Omega \rightarrow M$ be a measurable map. This function defines a probability measure $F_{*}P \in \mathcal{P}(M)$ by
	$$
	F_{*}P(A) := P(F^{-1}(A)) = P(\left\{\omega \in \Omega : F(\omega) \in A \right\}) \ \ \ \ \ (\forall A \in \mathcal{B}(M)), 
	$$ 
	which is called \textbf{pushforward measure} of $P$ by the function $F$. In probabilistic language, the pushforward measure $F_{*}P$ is called distribution of $F$. 
	
	\medskip
	
	Given $1\leq \theta\leq \infty$, we say that $F\in L^{\theta}(\Omega, M)$ if $F_*P\in \cP^\theta(M)$. In other words, for $1\leq \theta<\infty$, we say that $F\in L^{\theta}(\Omega, M)$ if for some (and hence for all) $y\in M$ it holds that
	\begin{equation}\label{cond}
	\int_\Omega \delta^\theta (F(\w), y)\, dP(\w)<\infty.
	\end{equation}
	On the other hand, we say that $F\in L^\infty(\Omega,M)$ if for some (and hence for all) $y\in M$ the function $\w\mapsto \delta(F(\w),y)$ is essentially bounded. 
	%

	\subsection{The inductive mean}\label{inductive}
	
	As in the case of strictly positive matrices considered in the introduction, we define the inductive means in general Hadamard spaces as follows.
	
	\begin{fed}(Inductive mean). Let $(M, \delta)$ be a Hadamard space. Given $a \in M^{\N}$ set
		\begin{align*}
		S_1(a) & = a_1  \\
		S_{n}(a) & = S_{n-1}(a)  \#_{\frac{1}{n}} a_{n} \ \ \ \ \ (n \geq 2).
		\end{align*}
	\end{fed}
	
	\begin{exa}\label{example inductivo}
		Suppose that $M$ is $\C^n$ with the usual euclidean distance. Since the geodesics in this case are the line segments, if we take a sequence $\{a_n\}_{n\in\N}$ in $\C^n$, then 
		\begin{align*}
		S_1(a) & = a_1,\\
		S_2(a) & = \frac{a_1+a_2}{2},\\
		S_3(a)&=\frac{2}{3}\left(\frac{a_1+a_2}{2}\right)+\frac13\, a_3=\frac{a_1+a_2+a_3}{3},
		\end{align*}
		and so on and so forth. Therefore, in this case the inductive means coincides the the arithmetic means.\EOE
	\end{exa}
	\bigskip
	
	From now on, let $(M,\delta)$ be a Hadamard space. As a consequence of \eqref{eq5}, we directly get the following result. 
	
	\begin{cor}\label{conse}
		Given $a,b \in M^\N$, then
		\begin{equation}\label{desig}
		\delta(S_n(a), S_n(b)) \leq \frac{1}{n}\sum_{i=1}^{n} \delta(a_i, b_i).
		\end{equation}
	\end{cor}
	
	\medskip
	
	The next lemma follows from \eqref{eq2}, and it is a special case of a weighted inequality considered by Lim and P\'alfia in \cite{[Lim2]}.
	
	\begin{lem}\label{lemma1}
		Given $a \in M^\N$ and $z\in M$, for every $k,m \in \N$ 
		\begin{align*}
		\delta^2(S_{k+m}(a), z)  &\leq \ \frac{k}{k+m}\ \delta^2(S_{k}(a), z) + \ds\frac{1}{k+m}\ds\sum_{j = 0}^{m - 1}\delta^2(a_{k+j+1}, z)\\ 
		&\quad - \ds\frac{k}{(k+m)^2}\ds\sum_{j = 0}^{m - 1}\delta^2(S_{k+j}(a), a_{k+j+1}). 
		\end{align*}
	\end{lem}
	
	\bdem
	By the inequality \eqref{eq2} applied to $S_{n+1}(a)=S_n(a)\,\#_{n+1}\,(a_{n+1})$ we obtain
	\begin{align*}
	(n+1)\ \delta^2(S_{n+1}(a), z) - n\ \delta^2(S_{n}(a), z) & \leq  \delta^2(a_{n+1}, z) - \ds\frac{n}{(n + 1)}\delta^2(S_{n}(a), a_{n+1}).
	\end{align*}
	Summing these inequalities from $n=k$ until $n=k+m-1$ we get that the difference 
	\begin{align*}
	(k+m)\ \delta^2(S_{k+m}(a), z) - k\ \delta^2(S_{k}(a), z), 
	\end{align*}
	obtained from the telescopic sum of the left hand side, is less or equal than 
	\begin{align*}
	\sum_{j=0}^{m-1}\left(\delta^2(a_{k+j+1}, z) - \ds\frac{k+j}{(k+j+1)}\delta^2(S_{k+j}(a), a_{k+j+1})\right).
	\end{align*}
	Finally, using that $\frac{k+j}{k+j+1}\geq \frac{k}{k+m}$ for every $j\in\{0,\ldots,m-1\}$, this sum is bounded from the above by 
	\begin{align*}
	\sum_{j=0}^{m-1}\left(\delta^2(a_{k+j+1}, z) - \ds\frac{k}{(k+m)}\delta^2(S_{k+j}(a), a_{k+j+1})\right),
	\end{align*}
	which completes the proof.
	\edem
	
	\medskip
	
	\noindent Given a sequence $a \in M^\N$, let $\Delta(a)$ denote the diameter of its image, that is
	$$
	\Delta(a) := \sup_{n,m\in\N} \delta(a_n, a_m).
	$$
	Note that, also by \eqref{eq2}, $\delta(S_n(a), a_k) \leq \Delta(a)$ for all $n, k \in \N$. 
	
	\medskip
	
	\begin{lem}\label{lemma2}
		Given $a\in M^\N$ such that $\Delta(a) < \infty$, then for all $k,m \in \N$ it holds that 
		\begin{align*}
		\frac{1}{m} \sum_{j = 0}^{m-1} \delta^2(S_{k}(a), a_{k+j+1}) \leq   \tilde{R}_{m,k} + \frac{1}{m} \sum_{j = 0}^{m-1}\delta^2(S_{k+j}(a), a_{k+j+1}).
		\end{align*}
		where $\tilde{R}_{m,k}=\left(\ds\frac{m^2}{(k+1)^2} + 2\ds\frac{m}{k+1}\right) \Delta^2(a)$. 
	\end{lem}
	
	\bdem
	
	Note that by \eqref{desig} and all $k$,
	$$
	\delta(S_{k+j}(a), S_{k+j+1}(a)) \leq \frac{1}{k+j+1}\Delta(a).
	$$
	Hence
	
	\begin{align*}
	\delta(S_{k}(a), a_{k+j+1}) & \leq  \delta(S_{k}(a), S_{k+j}(a)) + \delta(S_{k+j}(a), a_{k+j+1}) \\
	& \leq \ds\sum_{h=1}^{j} \ds\frac{1}{k+h} \Delta(a) + \delta(S_{k+j}(a), a_{k+j+1}) \\
	& \leq  \ds\frac{m}{k+1} \Delta(a) + \delta(S_{k+j}(a), a_{k+j+1}).
	\end{align*}
	Therefore, , for every $j\leq m$
	\begin{align*}
	\delta^2(S_{k}(a), a_{k+j+1}) & 
	\leq  \left(\ds\frac{m^2}{(k+1)^2} + 2\ds\frac{m}{k+1}\right) \Delta^2(a) + \delta^2(S_{k+j}(a), a_{k+j+1}),
	\end{align*}
	where we have used that $\delta(S_{k+j}(a), a_{k+j+1}) \leq \Delta(a)$ for every $k,j\in\N$. Summing up these inequalities and dividing by $m$, we get the desired result.
	\edem

	\section{Ergodic formulae associated to inductive means}
	
	\subsection{The framework and basic notation}
	
	Let $(G,+)$ be a compact, abelian, and metrizable topological group. In this group we fix a Haar measure $m$, and we take an ergodic automorphism $\tau(h)=h+g$ for some $g\in G$. A shift invariant metric $G$ is denoted by ${d}_{G}$. 
	Throughout this section we work with the dynamical system $(G,\tau)$.
	
	\medskip
	
	\begin{rem}
		A topological dynamical system $(\W,\tau)$ is called a \textbf{Kronecker system} if it is isomorphic to a group dynamical system $(G,\tau)$ as the one described above. Every topological Kronecker system $(\Omega; x \rightarrow x + \alpha)$ can be canonically converted into a measure-preserving system which is compact. It is well known that any isometric (or equicontinuous) and minimal dynamical system is a Kronecker system. 
	\end{rem}
	
	On the other hand, we will fix a Hadamard space $(M,\delta)$. Given a function $A: G \rightarrow M$, we define $a^{\tau} : G \rightarrow M^\N$ by
	\begin{equation}\label{eq def de suc}
	a^{\tau}(x) :=\{a^{\tau}_j(x)\}_{n\in\N}  \peso{where} a^{\tau}_j(x)=A(\tau^j(x)).
	\end{equation}
	
	\subsection{The continuous case}
	
	In this section we will prove the ergodic formula for continuous functions.
	
	\begin{teo}\label{mainteo1}
		Let $M$ be a Hadarmard space and $A:G\to M$ a continuous function. Then
		\begin{equation}
		\lim_{n\to\infty} S_{n}(a^{\tau}(g)) = \barza{A},
		\end{equation}
		uniformly in $g\in G$.
	\end{teo}
	
	With this aim, we firstly prove some technical results, which at the end of this subsection are combined to get a proof of Theorem \ref{mainteo1}.
	
	\begin{lem}\label{equicont}
		
		Let $A: G\rightarrow M$ be a continuous function, and let $K$ be any compact subset of $M$. For each $n\in\N$, define $F_n :G \times K \rightarrow \R$ by
		$$
		F_n(g, x) = \frac{1}{n} \sum_{j = 0}^{n-1} \delta^{2}(a_j^{\tau}(g), x).
		$$
		Then, the family $\left\{F_n\right\}_{n \in \N}$ is equicontinuous.
		
	\end{lem}

	\bdem
	

	By the triangular inequality, the map $y\mapsto \delta^2(A(\cdot),y)$ is continuous from $(K,\delta)$ into the set of real valued continuous functions defined on $G$ endowed with the uniform norm. Since $K$ is compact, the family $\{\delta^2(A(\cdot),x)\}_{x\in K}$ is (uniformly) equicontinuous. Hence, given $\eps>0$, there exists $\delta>0$ such that if
	$$
	d_G(g_1,g_2) < \delta \peso{then} |\delta^2(A(g_1),x)-\delta^2(A(g_2),x)|<\frac{\eps}{2},
	$$
	for every $x\in K$. Since $\tau$ is isometric and $d_G(g_1,g_2) < \delta$, we get that
	$$
	|F_n(g_1, x)-F_n(g_2, x)| = \left|\frac{1}{n} \sum_{j = 0}^{n-1} \delta^{2}(a_j^{\tau}(g_1), x)-\delta^{2}(a_j^{\tau}(g_2), x)\right|<\frac{\eps}{2}.
	$$
	Let $\Delta$ be the diameter of the set $\,(\mbox{Image}(A)\times K)$ in $M^2$. Since both sets are compact, $\Delta<\infty$. So, take $(g_1,x_1)$ and $(g_2,x_2)$ such that $d_G(g_1,g_2)<\delta$ and $\delta(x_1,x_2)<\frac{\eps}{4\Delta}$. Then
	\begin{align*}
	|F_n(g_1, x_1)-F_n(g_2, x_2)|&\leq  |F_n(g_1, x_1)-F_n(g_1, x_2)|+|F_n(g_1, x_2)-F_n(g_2, x_2)|\\
	&\leq \frac{2\Delta}{n}\sum_{k=0}^{n-1} \delta(x_1,x_2) \ +\ \frac{\eps}{2}<\eps.
	\end{align*}
	\edem
	
	Now, as a consequence of Arzel\`a-Ascoli and Birkhoff theorems we get
	
	\begin{pro}\label{pro10}
		Let $A:G\to M$ be a continuous function, and $K$ a compact subset of $M$. Then
		$$
		\lim_{n \rightarrow \infty}  \frac{1}{n} \sum_{j = 0}^{n-1} \delta^{2}(a_j^{\tau}(g), x) = \int_{G} \delta^2(A(\gamma), x)dm(\gamma),
		$$
		and the convergence is uniform in $(g, x) \in G \times K$.
	\end{pro}

	\bigskip
	
	From now on we will fix the continuous function $A:G\to M$. Let
	$$
	\alpha := \min_{x\in M} \int_{G} \delta^2(A(g), x)dm(g),
	$$
	and $\barza{A}$ is the point where this minimum is attained, i.e., $\barza{A}$ is the barycenter of the pushforward by $A$ of the Haar measure in $G$. 
	Then we obtain the following upper estimate.
	
	\medskip
	
	\begin{lem}\label{lemerg1}
		For every $\varepsilon > 0$, there exists $m_0 \in \N$ such that, for all $m \geq m_0$ and for all $k \in \N$,
		\begin{align*}
		\delta^2(S_{k+m}(a^{\tau}(g)), \barza{A})  & \leq \frac{k}{k+m}\delta^2(S_{k}(a^{\tau}(g)), \barza{A}) + \ds\frac{m}{k+m}(\alpha+\eps)  \\ 
		&\quad - \ds\frac{km}{(k+m)^2}\left(\frac{1}{m}\ds\sum_{j = 0}^{m - 1}\delta^2\big(S_{k+j}(a^{\tau}(g)), a^{\tau}_{k+j+1}(g)\big)\right). 
		\end{align*}
	\end{lem}
	
	\bdem
	For every $\varepsilon > 0$, exists $m_0 \in \N$ such that for all $m \geq m_0$,
	$$
	\left|\frac{1}{m}\ds\sum_{j = 0}^{m - 1}\delta^2(a_{k+j+1}^{\tau}(g), \barza{A}) - \alpha\right| < \varepsilon.
	$$
	Note that $m_0$ is independient of $k$ by Proposition \ref{pro10}. Now, by Lemma \ref{lemma1}
	\begin{align*}
	\delta^2(S_{k+m}(a^{\tau}(g&)),\barza{A})   \leq \frac{k}{k+m}\ \delta^2(S_{k}(a^{\tau}(g)), \barza{A}) + \ds\frac{1}{k+m}\ds\sum_{j = 0}^{m - 1}\delta^2(a_{k+j+1}^{\tau}(g), \barza{A})\\ 
	& \quad - \ds\frac{k}{(k+m)^2}\ds\sum_{j = 0}^{m - 1}\delta^2(S_{k+j}(a^{\tau}(g)), a_{k+j+1}^{\tau}(g)) \\
	& = \frac{k}{k+m} \ \delta^2(S_{k}(a^{\tau}(g)), \barza{A}) + \ds\frac{m}{k+m}\left(\frac{1}{m}\ds\sum_{j = 0}^{m - 1}\delta^2(a_{k+j+1}^{\tau}(g), \barza{A})\right) \\ 
	& \quad - \ds\frac{km}{(k+m)^2}\left(\frac{1}{m}\ds\sum_{j = 0}^{m - 1}\delta^2(S_{k+j}(a^{\tau}(g)), a_{k+j+1}^{\tau}(g))\right) \\
	& \leq \frac{k}{k+m}\ \delta^2(S_{k}(a^{\tau}(g)), \barza{A}) + \ds\frac{m}{k+m}(\alpha+\varepsilon)  \\ 
	& \quad - \ds\frac{km}{(k+m)^2}\left(\frac{1}{m}\ds\sum_{j = 0}^{m - 1}\delta^2(S_{k+j}(a^{\tau}(g)), a_{k+j+1}^{\tau}(g))\right).
	\end{align*}
	
	\edem
	
	\noindent Recall that, given a sequence $a \in M^\N$, then $\Delta(a)$ denotes the diameter of its image, i.e.,
	$$
	\Delta(a) := \sup_{n,m\in\N} \delta(a_n, a_m).
	$$
	Since $A:G\to M$ is continuous, note that  
	$$
	C_a:= \sup_{g\in G} \Delta(a^{\tau}(g))<\infty.
	$$

	\begin{lem}\label{lemerg2}
		
		For every $\varepsilon > 0$, there exists $m_0 \in \N$ such that for all $m \geq m_0$ and for all $k \in \N$
		$$
		\delta^2(S_k(a^{\tau}(g)), \barza{A}) - \varepsilon + \alpha - R_{m,k} \leq \frac{1}{m} \sum_{j = 0}^{m-1}\delta^2(S_{k+j}(a^{\tau}(g)), a_{k+j+1}^{\tau}(g)),
		$$
		where $\ds R_{m,k}= \left(\ds\frac{m^2}{(k+1)^2} + 2\ds\frac{m}{k+1}\right) C_a^2$.
	\end{lem}
	
	\bdem	
	Consider the compact set
	$$
	K := \overline{cc\left\{S_k(a^{\tau}(x)) : k \in \N\right\}},
	$$
	where the convex hull is in the geodesic sense.  For every $\varepsilon > 0$, there exists $m_0 \in \N$ such that for all $m \geq m_0$, by the variance inequality (Proposition~\ref{baryvar}) and Proposition \ref{pro10}, it holds that
	\begin{align*}
	\delta^2(S_k(a^{\tau}(x)), \barza{A}) & \leq \int_G \delta^2\big(S_k(a^{\tau}(g)), A(\gamma)\big)\,dm(\gamma)\ -\ \alpha \\
	&\leq \eps+ \frac{1}{m} \sum_{j = 0}^{m-1} \delta^2(S_{k}(a^{\tau}(g)), a_{k+j+1}^{\tau}(g)) - \alpha.
	\intertext{Finally, by Lemma \ref{lemma2}}
	\delta^2(S_k(a^{\tau}(x)), \barza{A}) & \leq \varepsilon  + \frac{1}{m} \sum_{j = 0}^{m-1}\delta^2(S_{k+j}(a^{\tau}(x)), a_{k+j+1}^{\tau}(x))+R_{m,k}- \alpha, 
	\end{align*}
	where $\ds R_{m,k}=\left(\frac{m^2}{(k+1)^2} + 2\ds\frac{m}{k+1}\right) C_a$.
	\edem

	\begin{lem}\label{lemerg3}
		Given $\eps > 0$, there exists $m_0\geq 1$ such that for every $\ell\in\N$
		$$
		\delta^2(S_{\ell m_0}(a^{\tau}(g)), \barza{A}) \leq \frac{L}{\ell} + \eps,
		$$
		uniformly in $g\in G$, where $L = \alpha + 3  C_a^2$.  
	\end{lem}
	
	\medskip
	
	\bdem
	Fix $\eps > 0$. By Lemmas \ref{lemerg1} and \ref{lemerg2}, there exists $m_0 \geq 1$ such that for all $k \in \N$,
	\begin{align*}
	\delta^2(S_{k+m_0}(a^{\tau}(g)), \barza{A})  
	& \leq \frac{k}{k+m_0}\delta^2(S_{k}(a^{\tau}(g)), \barza{A}) + \ds\frac{m_0}{k+m_0}\big(\alpha+\eps\big) \\ 
	&\quad - \ds\frac{km_0}{(k+m_0)^2}\left(\frac{1}{m_0}\ds\sum_{j = 0}^{m_0 - 1}\delta^2(S_{k+j}(a^{\tau}(x)), a_{k+j+1}^{\tau}(x))\right), 
	\end{align*}
	and
	\begin{align*}
	\frac{1}{m_0} \sum_{j = 0}^{m_0-1}\delta^2(S_{k+j}(a^{\tau}(g)), a_{k+j+1}^{\tau}(g)) \geq \delta^2(S_k(a^{\tau}(g)), \barza{A}) - \eps+ \alpha -  R_{m_0,k}.
	\end{align*}
	Therefore, combining these two inequalities we obtain
	\begin{align*}
	\delta^2(S_{k+m_0}(a^{\tau}(g)), \barza{A})  
	& \leq \frac{k}{k+m_0}\delta^2(S_{k}(a^{\tau}(g)), \barza{A}) + \ds\frac{m_0}{k+m_0}\big(\alpha+\eps\big) \\ 
	& \quad - \ds\frac{km_0}{(k+m_0)^2}\left(\delta^2(S_k(a^{\tau}(g)), \barza{A}) - \eps + \alpha -  R_{m_0,k}\right). 
	\end{align*}
	Consider now the particular case where $k=\ell m_0$. Since $\ds R_{m_0,\ell m_0} \leq \frac{3}{\ell} C_a^2$ we get
	
	\begin{align}
	\delta^2(S_{(\ell+1)m_0}(a^{\tau}(g)), \barza{A})  
	& \leq \frac{\ell}{\ell+1}\delta^2(S_{\ell m_0}(a^{\tau}(g)), \barza{A}) + \ds\frac{1}{l+1}\big(\alpha+\eps\big)-\nonumber \\ 
	& \quad \frac{\ell}{(\ell+1)^2}\left(\delta^2(S_{\ell m_0}(a^{\tau}(g)), \barza{A}) - \eps + \alpha -  R_{m_0,\ell m_0}\right)\nonumber\\
	&\leq \frac{\ell^2\ \delta^2(S_{\ell m_0}(a^{\tau}(g)), \barza{A}) +(2\ell+1)\eps+\alpha+3C_a^2}{(\ell+1)^2}.\label{eq recursiva}
	\end{align}
	Using this recursive inequality, the result follows by induction on $\ell$. Indeed, if $\ell=1$ then 
	$$
	\delta^2(S_{m_0}(a^{\tau}(g)),\barza{A})\leq C_a^2 \leq L.
	$$
	On the other hand,  if we assume that the result holds for some $\ell\geq 1$, i.e.
	$$
	\delta^2(S_{\ell m_0}(a^{\tau}(x)), g) \leq \frac{L}{\ell} + \eps,
	$$
	then combining this inequality with \eqref{eq recursiva} we have that
	\begin{align*}
	\delta^2(S_{(\ell+1)m_0}(a^{\tau}(g)), \barza{A})  
	&\leq \frac{\ell L+ \ell^2\eps +(2\ell+1)\eps+\alpha+3C_a^2}{(\ell+1)^2}=\frac{L}{\ell+1}+\eps.
	\end{align*}
	\edem
	
	\medskip
	
	\noindent Now we are ready to prove the ergodic formula for continuous functions.
	
	\bdem[of Theorem \ref{mainteo1}]
	Given $\eps > 0$, by Lemma \ref{lemerg3}, there exists $m_0 \in \N$ such that,
	$$
	\delta^2(S_{\ell m_0}(a^{\tau}(g)), \barza{A}) \leq \frac{L}{\ell} + \frac{\eps^2}{8},
	$$
	for every $\ell\in \N$ . Take $\ell_0 \in \N$ such that for all $\ell \geq \ell_0$,
	\begin{equation}\label{eq demo continua}
	\delta^2(S_{\ell m_0}(a^{\tau}(g)), \barza{A}) \leq  \frac{\eps^2}{4}.
	\end{equation}
	Let $n=\ell m_0+d$ such that $\ell\geq \ell_0$ and $d\in\{1,\ldots, m_0-1\}$. Since $x \#_{t} x = x$ for all $x \in M$, using Corollary \ref{conse} with the sequences
	\begin{align*}
	&(\ a^\tau_1(g),\ldots,a^\tau_{\ell m_0}(g), \underbrace{S_{\ell m_0}(a^{\tau}(g)),\ldots, S_{\ell m_0}(a^{\tau}(g))}_{\mbox{\tiny{d times}}}\ )\\
	\intertext{and}
	&(\ a^\tau_1(g),\ldots,a^\tau_{\ell m_0}(g), \ a^\tau_{\ell m_0+1}(g)\ ,\ \ldots\ ,\ a^\tau_{\ell m_0+d}(g)\ ).
	\end{align*}
	Taking into account that $\delta(S_{\ell m_0}(a^{\tau}(g)),a^\tau_{\ell m_0+j}(g))\leq C_a$ for every $j\in\{1,\ldots,m_0-1\}$, we get
	\begin{align*}
	\delta(S_{\ell m_0}(a^{\tau}(g)), S_{\ell m_0+d}(a^{\tau}(g))) 
	& \leq \frac{1}{\ell m_0 + d}\sum_{j = 1}^{d} \delta(S_{\ell m_0}(a^{\tau}(g)),a^\tau_{\ell m_0+j}(g))\\
	& \leq \frac{d}{\ell m_0 + d} C_a 
	\leq\frac{1}{\ell} C_a\xrightarrow[k \rightarrow \infty]{} 0.
	\end{align*}
	Combining this with \eqref{eq demo continua} we obtain that for $n$ big enough $\delta(S_{n}(a^{\tau}(g)), \barza{A})<\eps$.
	\edem

	\subsection{The $L^{1}$ case}
	
	The natural framework for the ergodic theorem is $L^1$. In this section we will prove the following ergodic theorem for functions in $L^1(G,M)$ in terms of the inductive means, which is the main result of this paper. 
	
	\begin{teo}\label{Ergodic L1}
		
		Given $A\in L^1(G, M)$, for almost every $g\in G$
		\begin{equation}\label{eq ergodic L1}
		\lim_{n\to\infty }S_{n}(a^{\tau}(g))=\barza{A}. 
		\end{equation}
	\end{teo}
	
	\medskip
	
	\begin{rem}
		Let $M$ be the field of complex numbers with the usual distance. As we observed in the Example \ref {example inductivo}, the inductive means $S_{n}(a^{\tau}(g))$ become the usual arithmetic mean
		$$
		\frac{1}{n}\sum_{k=0}^{n-1} A(\tau^{n}(g)).
		$$
		On the other hand, the barycenter $\barza{A}$ is just the integral of $A$ with respect to the Haar measure $m$. So, equation \eqref{eq ergodic L1} takes the form
		$$
		\lim_{n\to\infty }\frac{1}{n}\sum_{k=0}^{n-1} A(\tau^{n}(g))=\int_G A(g)\, dm(g)
		$$ 
		which is the usual Birkhoff ergodic theorem. \EOE
	\end{rem}
	
	The strategy of the proof consists in constructing good approximations by continuous functions, and get the result of $L^1$ functions as a consequence of the theorem for continuous functions (Theorem \ref{mainteo1} above). So, the first questions that appear are: what does good approximation mean?, and what should we require to the approximation in order to get the $L^1$ case as a limit of the continuous case? The next two lemmas contain the clue to answer these two questions. 
	
	\bigskip

	\begin{lem}\label{lema 1 para L1}
		Let $(\W,\cB, P)$ be a probability space, and $A, B \in L^1(X, M)$. If
		\begin{align*}
		\barza{A} & = \argmin_{z \in M} \int_{\W} [\delta^2(A(\w), z) - \delta^2(A(\w), y)]\ dP(\w), \\
		\barza{B} & = \argmin_{z \in M} \int_{\W} [\delta^2(B(\w), z) - \delta^2(B(\w), y)]\ dP(\w),
		\end{align*}
		then 
		\begin{equation}\label{eq lema 1 L1}
		\delta(\barza{A}, \barza{B}) \leq  \int_\W \delta(A(\w), B(\w))dP(\w).
		\end{equation}
	\end{lem}
	
	\medskip
	
	\begin{rem}
		Recall that the definition of $\barza{A}$ (resp. $\barza{B}$) does not depend on the chosen $y \in M$.
	\end{rem}	
	
	\medskip
	
	\bdem
	
	By the variance inequality (Proposition \ref{baryvar}) we get
	\begin{align*}
	\delta^2(\barza{A}, \barza{B}) & \leq  \int_\W \delta^2(\barza{A}, B(\omega)) - \delta^2(\barza{B}, B(\omega))dP(\omega), \\
	\delta^2(\barza{A}, \barza{B}) & \leq \int_\W \delta^2(\barza{B}, A(\omega)) - \delta^2(\barza{A}, A(\omega))dP(\omega),
	\end{align*}
	and the combination of  these two inequalities leads to
	\begin{align*}
	2\delta^2(\barza{A}, \barza{B}) & \leq \int_\W \delta^2(\barza{A}, B(\omega)) + \delta^2(\barza{B}, A(\omega))\\
	& \quad - \delta^2(\barza{B}, B(\omega))  - \delta^2(\barza{A}, A(\omega))dP(\omega).
	\end{align*}
	Finally, using the Reshetnyak's cuadruple comparison (Theorem \ref{Res}) we obtain 
	\begin{align*}
	2\delta^2(\barza{A}, \barza{B}) & \leq 2\delta(\barza{A}, \barza{B})\ \int_\W \delta(A(\omega), B(\omega))\,dP(\w),
	\end{align*}
	which is, after an algebraic simplification, the desired result. 
	\edem
	
	\bigskip
	
	\begin{lem}\label{lema 2 para L1}
		Let $A, B \in L^1(G, M)$.  Given $\eps>0$, for almost every $g\in G$ there exists $n_0$, which may depends on $g$, such that 
		\begin{equation}\label{eq lema 1 L1}
		\delta\big(S_{n}(a^{\tau}(g)), S_{n}(b^{\tau}(g))\big) \leq \eps + \int_{G} \delta(A(g), B(g))dm(g),
		\end{equation}
		provided $n\geq n_0$.
	\end{lem}
	
	\bdem
	Indeed, by Corollary \ref{conse}
	$$
	\delta\big(S_n(a^{\tau}(g)), S_n(b^{\tau}(g))\big) \leq \ds\frac{1}{n} \ds\sum_{k = 0}^{n-1} \delta\big(a_k^{\tau}(g), b_k^{\tau}(g)\big) 
	=  \ds\frac{1}{n} \ds\sum_{k = 0}^{n-1} \delta\big(A(\tau^k(g)), B(\tau^k(g))\big),
	$$
	and therefore, the lemma follows by Birkhoff ergodic Theorem.
	\edem

	\subsubsection{Good approximation by continuous functions}\label{aproximando}
	
	The previous two lemmas indicate that we need a kind of $L^1$ approximation. More precisely, given $A\in L^1(G,M)$ and $\eps>0$, we are looking for a continuous function $A_\eps:G\to M$ such that
	$$
	\int_G \delta(A(g),A_\eps(g))\,dm(g)<\eps.
	$$ 
	In some cases there exists an underlyng finite dimensional vector space. This is the case, for instance, when $M$ is the set of (strictly) positive matrices, or more generally, when $M$ is a Riemannian manifold with non-positive curvature. In these cases, the function $A_\eps$ 
	can be constructed by using mollifiers. This idea was used by Karcher in \cite{Karcher}. In the general case, we can use a similar idea. 
	
	\bigskip
	
	Given $\eta>0$, let $U_\eta$ be a neighborhood of the identity of $G$ so that $m(U_\eta)<\eta$, whose diameter is also less than $\eta$. Fix any $y\in M$, and define
	\begin{equation}\label{eq def las Aeps}
	A_\eta(g_0)=\argmin_{z\in M} \int_{U_\eta}  [\delta^2(z, A(g+g_0)) - \delta^2(y, A(g+g_0)) ]\ dm(g).
	\end{equation}
	
	Equivalently, $A_\eta(g_0)$ is the barycenter of the pushforward by $A$ of the Haar measure restricted to $g_0+U_\eta$. This definition follows the idea of mollifiers, replacing the arithmetic mean by the average induced by barycenters. We will prove that, as in the case of usual mollifiers, these continuous functions provide good approximation in $L^1$ (Theorem \ref{aproximacion por continuas} below). With this aim, firstly we will prove the following lemma.

	\begin{lem}\label{lema aprox}
		Let $A \in L^p(G, M)$ where $1 \leq p < \infty$. Define the function $\fii : G \rightarrow [0,+\infty)$  by 
		$$
		\varphi(h) = \int_{G} \delta^p\big(A(g), A(g + h)\big)dm(g),
		$$
		is a continuous function.	
	\end{lem}

	\bdem
	
	Fix $z_0 \in M$, and define the measure on the Borel sets of $G$
	$$
	\nu(B) := \int_{B} \delta^p(A(g), z_0)dm(g).
	$$
	By definition, $\nu$ is absolutely continuous with respect to the Haar measure $m$. In consequence, given $\varepsilon > 0$, 
	there exists $\eta > 0$, such that, whenever a Borel set $B$ satisfies
	$$
	\int_{B} dm(g) < \eta,
	$$
	it holds that
	\begin{equation}\label{medida}
	\nu(B) =  \int_{B} \delta^p(A(g), z_0)dm(g) < \frac{\varepsilon}{2^{p+2}},
	\end{equation}
	By Lusin Theorem \cite[Thm 7.5.2]{Dudley}, there is a compact set $C_{\eta} \subset G$ such that $m(C_{\eta}) \geq 1 - \eta/2$ and the restriction of $A$ to $C_{\eta}$ is  (uniformly) continuous. 
	
	\medskip
	
	Since $m$ is a Haar measure, it is enough to prove the continuity of $\fii$ at the identity. With this aim, take a neighborhood of the identity $U$ so that whenever $g_1,g_2\in C_\eta$ satisfy that $g_1-g_2\in U$, it holds that
	$$
	\delta^p(A(g_1), A(g_2)) \leq \frac{\varepsilon}{2}.
	$$ 
	Given $h\in U$, define $\W:= C_{\eta} \cap (C_{\eta} + h)$, and $\W^c :=G \setminus \W$. Then
	\begin{align*}
	\int_G \delta^p(A(g), A(g+ h))dm(g) & = \int_{\Omega}+\int_{\Omega^c} \delta^p(A(g), A(g + h))dm(g) 
	\\ 
	&  \leq \frac{\varepsilon}{2}+ \int_{\Omega^c} \delta^p(A(g), A(g+h))dm(g) \\
	& \leq \frac{\varepsilon}{2} +  \int_{\Omega^c} [\delta(A(g), z_0) + \delta(A(g + h), z_0)]^p\, dm(\omega) \\
	& = \frac{\varepsilon}{2} +  2^{p+1} \int_{\Omega^c} \delta^p(A(g), z_0)dm(g),
	\end{align*}
	where in the last identity we have used that $m$ is shift invariant. Since $|\W^c|<\delta$ we obtain that 
	\begin{align*}
	\int_G \delta^p(A(g), A(g+ h))dm(g) & \leq  \frac{\varepsilon}{2}+ \frac{\varepsilon}{2}=\varepsilon.
	\end{align*}
	\edem
	
	\begin{cor}\label{son continuas}
		For every $\eta>0$, the functions $A_\eta$ are continuous.
	\end{cor}
	\bdem
	Indeed, by Lemma \ref{lema 1 para L1}
	\begin{align*}
	\delta(A_\eta(h_1),A_\eta(h_2))&\leq\frac{1}{m(U_\eta)}\int_{U_\eta} \delta(A(g+h_1),A(g+h_2))\,dm(g)\\
	&\leq\frac{1}{m(U_\eta)}\int_G \delta(A(g+h_1),A(g+h_2))\,dm(g)\\
	&\leq\frac{1}{m(U_\eta)}\int_G \delta[\,A(g),A(g+(h_2-h_1))\,]\,dm(g).
	\end{align*}
	So, the continuity of $A_\eta$ is a consequence of the continuity of $\fii$ at the identity.
	\edem
	
	\medskip
	
	\noindent The map $A\mapsto A_\eps$ has the following useful continuity property.
	
	\begin{lem}\label{supremo}
		Let $A, B \in L^1(G, M)$, and $\eta > 0$. For every $\eps>0$, there exists $\rho>0$ such that if
		$$
		\int_G \delta(A(g), B(g))dm(g) \leq \rho,
		$$
		then the corresponding continuous functions $A_\eta$ and $B_\eta$ satisfy that
		$$
		\max_{g \in G} \, \delta(A_{\eta}(g), B_{\eta}(g)) \leq \eps.
		$$
	\end{lem}
	\bdem
	Indeed, given $\eps>0$, take $\rho=m(U_{\eta})\eps$. Then, by Lemma \ref{lema 1 para L1}
	\begin{align*}
	\delta(A_{\eta}(g), B_{\eta}(g)) & \leq  \frac{1}{\left| U_{\eta} \right| }\int_{U_{\eta}} \delta(A(g + h), B(g + h))dm(h)\\
	& \leq \frac{1}{\left| U_{\eta} \right| }\int_{G} \delta(A(h), B(h))dm(h) \leq \eps,
	\end{align*}
	for all $g \in G$. 
	\edem
	
	We arrive to the main result of this section. 
	
	\begin{pro}\label{aproximacion por continuas}
		Given a function $A\in L^1(G,M)$, if $A_\eta$ are the continuous functions defined by \eqref{eq def las Aeps} then
		$$
		\lim_{\eta\to0^+}\int_G\delta(A(g),A_\eta(g))\,dm(g) = 0.
		$$
	\end{pro}
	\bdem
	Firstly, assume that $A\in L^2(G,M)$. In this case, by the variance inequality,  it holds that
	$$
	\delta^2(A(g), A_{\eta}(g)) \leq \frac{1}{\left| U_{\eta} \right| }\int_{U_{\eta}} \delta^2(A(g), A(g + h))dm(h).
	$$
	So, using Fubini's theorem we get
	\begin{align*}
	\int_G\delta^2(A(g),A_\eta(g))\,dm(g) &\leq \frac{1}{\left| U_{\eta} \right| }\int_{U_{\eta}} \int_G \delta^2(A(g), A(g + h))\, dm(g)\, dm(h)\\
	&= \frac{1}{\left| U_{\eta} \right| }\int_{U_{\eta}} \fii(h) \ dm(h).
	\end{align*}
	By Lemma \ref{lema aprox}, the function $\fii$ is continuous. In consequence, if $e$ denotes the identity of $G$
	$$
	\lim_{\eta\to 0^+} \frac{1}{\left| U_{\eta} \right| }\int_{U_{\eta}} \fii(h) \ dm(h) = \fii(e)=0.
	$$
	This proves the result for functions in $L^2(G,M)$ since by Jensen's inequality
	$$
	\ \ \int_G\delta(A(g),A_\eta(g))\,dm(g)\leq \left( \int_G\delta^2(A(g),A_\eta(g))\,dm(g) \right)^{1/2}.
	$$
	Now, consider a general $A \in L^1(G, M)$. Fix $z_0\in M$, and for each natural number $N$ define the truncations
	$$
	A^{(N)}(g) := \begin{cases}
	A(g) & \mbox{if  $\delta(A(g), z_0) < N$} \\
	z_0 & \mbox{if  $\delta(A(g), z_0) \geq N$}
	\end{cases}.
	$$
	For each $N$ we have that $A^{(N)}\in L^1(G,M)\cap L^\infty (G,M)$, and therefore it also belongs to $L^2(G,M)$. 
	On the other hand, since the function defined on $G$ by $g\mapsto \delta(A(g),z_0)$ is integrable, it holds that
	\begin{equation}\label{eq truncaditos}
	\int_G \delta(A(g),A^{(N)}(g))\ dm(g)= \int_{\{g:\,\delta(A(g), z_0) \geq N\}} \delta(A(g),z_0)\ dm(g) \xrightarrow[N\to \infty]{}0
	\end{equation}
	So, if $A_\eta$ and $A^{(N)}_\eta$ are the continuous functions associated to $A$ and $A^{(N)}$ res\-pec\-tively, then
	\begin{align*}
	\int_G \delta(A(g),A_\eta(g))\ dm(g)  &\leq \int_G \delta(A(g),A^{(N)}(g))\ dm(g) \\ 
	&\quad + \int_G \delta(A^{(N)}(g),A^{(N)}_\eta(g))\ dm(g) \\
	&\quad + \int_G \delta(A^{(N)}_\eta(g),A_\eta(g))\ dm(g). 
	\end{align*}
	Note that each term of the right hand side tends to zero: the first one by \eqref{eq truncaditos}, the second one by the $L^2$ case done in the first part, and the last one by Lemma \ref{supremo}. 
	\edem
	
	\subsubsection{Proof of Theorem \ref{Ergodic L1}}\label{se termina}
	
	Let $\eps>0$.  For each $k\in \N$ let $A_k$ be a continuous function such that
	$$
	\int_G \delta(A(g),A_k(g))\,dm(g)\leq\frac{1}{k}.
	$$
	By Lemma \ref{lema 2 para L1}, we can take a set of measure zero $N\subseteq G$ such that if we take $g\in G\setminus N$ and $k\in \N$, 
	there exists $n_0$, which may depend on $g$ and $k$, so that 
	$$
	\delta\big(S_{n}(a^{\tau}(g)), S_{n}(a^{\tau}_{(k)}(g))\big) \leq \frac{\eps}{4} + \int_{G} \delta(A(g), A_k(g))dm(g),
	$$
	provided $n\geq n_0$. In this expression, $a_{(k)}^{\tau}$ is the sequence defined in terms of $A_k$ and $\tau$ as in \eqref{eq def de suc}. Fix $g\in G\setminus N$. Taking $k$ so that $1/k<\eps/4$, we get that 
	$$
	\delta\big(S_{n}(a^{\tau}(g)), S_{n}(a^{\tau}_{(k)}(g))\big) \leq \frac{\eps}{2},
	$$
	for every $n\geq n_0$. By Lemma \ref{lema 1 para L1}, it also holds that
	$
	\delta(\barza{A}, \barza{A_k}) \leq  \frac{\eps}{4}
	$
	where
	\begin{align*}
	\barza{A} & = \argmin_{z \in M} \int_G [\delta^2(A(g), z) - \delta^2(A(g), y)]\ dm(g), \\
	\barza{A_k} & = \argmin_{z \in M} \int_G [\delta^2(A_k(g), z) - \delta^2(A_k(g), y)]\ dm(g).
	\end{align*}
	
	Finally, by Theorem \ref{mainteo1}, there exists $n_1\geq 1$ such that for every $n\geq n_1$
	$$
	\delta(S_{n}(a_{(k)}^{\tau}(g),\barza{A_k})\leq\frac{\eps}{4}.
	$$
	Combining all these inequalities we obtain that
	\begin{align*}
	\delta\big(S_{n}(a^{\tau}(g)), \barza{A}\big)&\leq \delta\big(S_{n}(a^{\tau}(g)), S_{n}(a^{\tau}_{(k)}(g))\big)\\
	&\quad +\delta\big(S_{n}(a^{\tau}_{(k)}(g)),\barza{A_k}\big)+\delta(\barza{A_k}, \barza{A})\leq\eps,
	\end{align*}
	which concludes the proof.
	
	\section*{Acknowledgements:}
	This work was supported by Consejo Nacional de Investigaciones Cient\'\i ficas y T\'ecnicas-Argentina (PIP-152), Agencia Nacional de Promoci\'on de Ciencia y Tecnolog\'\i a-Argentina (PICT 2015-1505), Universidad Nacional de La Plata-Argentina (UNLP-11X585) and  Ministerio de Economía y Competitividad-Espa\~na (MTM2016-75196-P).


\begin{thebibliography}{XX}
		
		\bibitem{Alex}  Alexandrov, A. D.: A theorem on triangles in a metric space and some applications, Trudy Math. Inst. Steklov 38  (1951), 5-23.
		
		\bibitem{[Ando]}  Ando T.,  Li C. K.,  Mathias R.: Geometric means, Linear Algebra Appl. 385 (2004), 305-334.
		
		\bibitem{Austin}  Austin T., A CAT(0) valued pointwise ergodic theorem, J. Topol. Anal. 3 (2011) 145-152.
		
		\bibitem{ball} Ballmann W.: Lectures on spaces of nonpositive curvature. DMV Seminar Band 25, Birkh\"auser Verlag, Basel.
		
		\bibitem{BGS} Ballmann W.,  Gromov M.,  Schroeder V.: Manifolds of nonpositive curvature. Progress in Mathematics 61. Birkh\"auser Boston Inc., Boston, MA
		
		\bibitem{[Berger]} Berger M.: A Panoramic View of Riemannian Geometry, Springer, 2003.
		
		\bibitem{[Bhatia1]} Bhatia R.: Positive Definite Matrices, Princeton Series in Applied Mathematics, Princeton University Press, 2007.
		
		
		\bibitem{[Bhatia2]} Bhatia R., Karandikar R.: Monotonicity of the matrix geometric mean. Math. Ann. 353(4)(2012) 1453-1467.
		
		\bibitem{Bi-Ia}Bini D., Iannazzo B.: Computing the Karcher mean of symmetric positive definite matrices, Linear
		Algebra Appl. 438 (2013) 1700-1710.
		
		\bibitem{Bochi} Bochi J., Navas A.: A geometric path from zero Lyapunov exponents to rotation cocycles, Ergodic Theory Dynam. Systems 35 (2015) 374-402.
		
		
		\bibitem{BH} Bridson M.R. , Haefliger A.: Metric spaces of non-positive curvature, Grundlehren der Mathematischen Wissenschaften, 319. Springer-Verlag, Berlin.
		
		\bibitem{Barba} Barbaresco F.: Interactions between symmetric cone and information geometries: Bruhat-Tits and Siegel
		spaces models for higher resolution autoregressive Doppler imagery, Emerging Trends in Visual Computing,
		Lecture Notes in Computer Science 5416 (2009) 124-163.
		
		
		\bibitem{Dudley} Dudley R. M.:  Real analysis and probability, Cambridge Studies in Advanced Mathematics, 74. Cambridge University Press, Cambridge. 
		
		\bibitem{Sahib} Es-Sahib A.,  Heinich H.: Barycentre canonique pour un espace metrique ?a courbure negative, in: Seminaire de Probabilites, in: Lecture Notes in Math., vol. 1709, Springer, Berlin, 1999.
		
		
		\bibitem{[Holbrook]} Holbrook J.: No dice: a determinic approach to the Cartan centroid, J. Ramanujan Math. Soc. 27 (2012) 509-521.
		
		\bibitem{[Jost]} Jost J.: Nonpositive Curvature: Geometric and Analytic Aspects, Lectures in Mathematics ETH Zurich, Birkhauser, 1997.
		
		\bibitem{Karcher}  Karcher H.: Riemannian center of mass and mollifier smoothing, Comm. Pure Appl. Math. 30 (1977) 509-541.
		
		\bibitem{LL1}  Lawson J.,  Lim Y.:  A general framework for extending means to higher orders, Colloq. Math. 113 (2008) 191-221.
		
		\bibitem{[LyL]} Lawson J.,  Lim Y.: Monotonic properties of the least squares mean, Math. Ann. 351 (2011) 267-279.
		
		\bibitem{LL2} Lawson J.,  Lim Y.: Weighted means and Karcher equations of positive operators, Proc. Natl. Acad. Sci. USA 110 (2013) 15626-15632.
		
		\bibitem{LL3} Lawson J.,  Lim Y.: Karcher means and Karcher equations of positive definite operators, Trans. Amer. Math. Soc. Ser. B1 (2014) 1-22.
		
		\bibitem{LL4} Lawson J.,  Lim Y.: Contractive barycentric maps, J. Operator Theory 77 (2017) 87-107.
		
		\bibitem{LeeNaor} Lee J.,  Naor A.: Extending Lipschitz functions via random metric partitions, Invent. Math. 160 (2005) 59-95.
		
		\bibitem{[Lim]}  Lim Y.: Riemannian distances between Geometric means. SIAM J. Matrix Anal. Appl.  34 (2013), 932-945.
		
		\bibitem{LimPal}  Lim Y.,  P\'alfia M.: Matrix power mean and the Karcher mean, J. Funct. Anal. 262 (2012) 1498-1514.
		
		\bibitem{[Lim2]} Lim Y.,  P\'alfia M.: Weighted deterministic walks and no dice approach for the least squares mean on Hadamard spaces, Bull. Lond. Math. Soc. 46 (2014) 561-570.
		
		\bibitem{[Lim3]} Lim Y.,  P\'alfia M.: Approximations to the Karcher mean on Hadamard spaces via geometric power means, Forum Math. 27 (2015) 2609-2635.
		
		
		\bibitem{Mendel} Mendel  M., Naor  A.: Spectral calculus and Lipschitz extension for barycentric metric spaces, Anal. Geom. Metr. Spaces 1 (2013) 163-199.
		
		\bibitem{Maher}  Moakher M.,  Zerai M.: The Riemannian geometry of the space of positive-definite matrices and
		its application to the regularization of positive-definite matrix-valued data, J. Math. Imaging Vision 40
		(2011) 171-187.
		
		\bibitem{Navas}  Navas A.: An $L^1$ ergodic theorem with values in a non-positively curved space via a canonical barycenter map, Ergodic Theory Dynam. Systems 33 (2013) 609-623.
		
		\bibitem{Ohta}  Ohta S.: Extending Lipschitz and H\"older maps between metric spaces, Positivity 13 (2009) 407-425.
		
		\bibitem{[Palfia]}  P\'alfia M.: Means in metric spaces and the center of mass, J. Math. Anal. Appl. 381 (2011) 383-391.
		
		\bibitem{Pass1}  Pass B.: Uniqueness and Monge solutions in the multimarginal optimal transportation problem, SIAM J. Math. Anal. 43 (2011) 2758-2775.
		
		\bibitem{Pass2}  Pass B.: Optimal transportation with infinitely many marginals, J. Funct. Anal. 264 (2013) 947-963.
		
		\bibitem{Res}  Reshetnyak Y. G.: (1968): Inextensible mappings in a space of curvature no greater
		than K. Sib. Math. Jour. 9, 683-689
		
		\bibitem{[sturm]}  Sturm K.-T.: Probability measures on metric spaces of nonpositive curvature, in: P. Auscher, et al. (Eds.), Heat Kernels and Analysis on Manifolds, Graphs, and Metric Spaces, in: Contemp. Math., vol. 338, Amer. Math. Soc. (AMS), Providence,
		2003.
		
	\end{thebibliography}
\end{document}